\documentclass[a4paper,10pt]{article}

\usepackage{amsmath}
\usepackage{amsfonts}
\usepackage{amsthm}

\def\R{\mathbb{R}}
\def\C{\mathbb{C}}

\def\A{\mathcal{A}}

\def\Geh{\mathfrak{g}}

\def\U{\mathfrak{U}}
\def\S{\mathfrak{S}}
\def\Z{\mathfrak{Z}}

\def\lin{\text{Lin}}
\def\pol{\text{Pol}}

\def\deg{\text{deg}}
\def\CH{\text{\rm CH}}

\def\CS{\mathcal{S}}
\def\O{\mathcal{O}}

\def\os{\text{\rm Os}}
\def\osint{\text{\rm Os-}\int}

\newtheorem{defi}{Definition}[section]
\newtheorem{thm}{Theorem}[section]
\newtheorem{prop}{Proposition}[section]
\newtheorem{lem}{Lemma}[section]
\newtheorem{cor}{Corollary}[section]

\def\remark{{\noindent \em Remark:}}

\DeclareMathOperator{\Ad}{Ad}

\DeclareMathOperator{\SO}{SO}
\DeclareMathOperator{\so}{\mathfrak{so}}

\begin{document}

\title{Differentiability of quantum moment maps and 
$G$-invariant star products
}


\title{Differentiability of quantum moment maps and 
$G$-invariant star products
\footnote{ This research was supported by a Postdoctoral Scholarship of 
the Ministry of National Education and Research of France.}
}
\author{Kentaro Hamachi \\
\small \it 
Laboratoire Gevrey de Math\'ematique Physique, Universit\'e de Bourgogne, \\
\small \it 
 BP 47870, F-21078 Dijon Cedex, France. \\
\small \it hamachi@u-bourgogne.fr
}
\date{September 27, 2002}
\maketitle



\maketitle

\begin{abstract}
We study quantum moment maps of $G$-invariant star products, which are 
a quantum analogue of the moment map for classical Hamiltonian systems. 
Introducing an integral representation, we show that any quantum 
moment map for a $G$-invariant star product is differentiable. 
This property gives us a new method for the classification of 
$G$-invariant star products on regular coadjoint orbits of compact 
semisimple Lie groups.
\end{abstract}

\section{Introduction}
Deformation quantization was introduced by Bayen, Flato, 
Fronsdal, Lichnerowicz and Sternheimer in the 70's\ \cite{BFFLS:1977}. 
It constitutes one of the important methods for quantizing 
classical systems. 
This quantization scheme provides an autonomous theory based on 
deformations of the ring of classical observables on a phase space 
(Poisson algebra), and does not involve 
a radical change in the nature of the observables. 

Star products invariant under the action of a Lie group $G$ have 
been studied with increasing generality from the beginning 
of the deformation quantization. 
They appear naturally in the quantization of 
classical systems with group symmetries, or in the star 
representation theory of Lie groups. 

Quantum moment maps have been introduced in \cite{Xu:1998}, 
and are the natural quantum analogue of moment maps 
on Hamiltonian $G$-spaces \cite{JT:1994} ; 
see Definition~\ref{def:QMM}.                         
A quantum moment map plays an important role for the study of 
$G$-invariant star products, similar to the one played by 
(classical) moment map for classical systems. 
One of the interesting applications of quantum moment maps is to 
provide an example of {\em quantum dual pair} 
\cite{Xu:1998,Weinstein:1983}. 
Another remarkable result is the {\em quantum reduction theorem}, 
which says that a quantization commutes with reduction 
\cite{Fedosov:1998}. 
We also give an application of quantum moment map by providing 
an invariant, called $c_*$ in  \cite{Hamachi:1999},
for a $G$-invariant star product $*$ on a $G$-transitive 
symplectic manifold \cite{Hamachi:1999}. 
This $c_*$ is computed with the help of a quantum moment map and 
depends only on the class of $G$-equivalent star products. 
In \cite{Hamachi:1999,Hamachi:2002}, we give a few 
examples of $c_*$ for a $\SO(3)$-invariant star product on the 
coadjoint orbit $S^2$. 

But there are serious problems with quantum moment maps. 
First, there is no obvious way to compute an explicit expression 
for a quantum moment map for a given $G$-invariant star product. 
We provide a partial answer to this problem in \cite{Hamachi:2002}.
Another important problem is the differentiability of 
quantum moment maps. 
Originally, a quantum moment map is defined only on the universal 
enveloping algebra $\U(\Geh_\lambda)$, that is, the set of 
polynomials on $\Geh^*$. 
But this definition of quantum moment maps does not 
directly imply its differentiability. 
A priori, a quantum moment map has only an algebraic meaning, 
and cannot be studied in the category of differentiable 
deformations, which can be inconvenient. 

In this article, we give another expression for quantum moment maps 
which is differentiable. 
This expression is an analogue of Weyl correspondence that 
can be formally written as: 
\begin{align*}
\Phi_*(u)=\int {\mathfrak F}u(\xi)\exp_*(i\xi\Phi_*(X)) d\xi,
\end{align*}
where ${\mathfrak F}u$ denotes the Fourier transformation of $u$, 
and the $\Phi_*$ in the integral is a quantum moment map of $*$ on $\Geh$. 
To make sense of the above formula, it is necessary to address 
two questions, to define the function $\exp_*(i\xi\Phi_*(X))$ 
and to give a meaning to the integral. 

For the first one, we simply define $\exp_*(i\xi\Phi_*(X))$ by 
power series with respect to the star product. 
We show that this naive definition of $\exp_*(i\xi\Phi_*(X))$ 
is well defined and it is a product of 
$e^{i\xi\Phi_0(X)}$ and a polynomial in $\xi$. 
This is a ingredient to make the quantum moment map differentiable. 
For the second one, since the domain of a quantum moment map 
contains any polynomial, $u$ in the above formula should be 
considered as a tempered distribution. 
In fact, for any slowly increasing infinitely differentiable 
function $u$, one can provide the integration as 
${\mathfrak F}_x^{-1}[{\mathfrak F}_\xi[u](\xi)
\exp_*(i\xi\Phi_*(X))e^{-i\xi \Phi_0(X)}]|_{x=\Phi_0(X)}$ 
 \cite{Schwartz:1966}. 
We prefer to use oscillatory integrals rather than tempered 
distributions in order to make computation easier. 
We give a brief review on oscillatory integrals in Appendix; 
see also \cite{Kumano-go:1981}. 

As an application of the differentiability of a quantum moment map 
we give a structure theorem for $G$-invariant star products on 
a coadjoint orbit of compact semisimple Lie groups. 
The class of $G$-invariant star products is parametrized by 
$G$-invariant Weyl curvature, that is, 
the second $G$-invariant de Rham cohomology \cite{BBG:1998}. 
However this classification does not give enough information on the 
structure of these star products. 

In relation to the structure of star products, there is an 
interesting study in \cite{FL:2001}. 
It provides a family of algebraic star products on a coadjoint orbit 
of semisimple Lie group by a quotient algebra of the Gutt star 
product. 
This work has the advantage to give an explicit representation 
of such kind of star products. 

We provide here a similar structure theory of $G$-invariant 
star products on such orbits in the differentiable category as an 
application of the differentiability of quantum moment maps. 
So we have another classification of such star products by using 
quantum moment maps. 
Moreover, as a corollary of the structure theorem, we have an 
answer of the problem we introduced in \cite{Hamachi:2002}: 
``Does $c_*$ parameterize the class of $G$-invariant star products?''
The answer is ``yes'' for regular coadjoint orbits of compact semisimple 
Lie groups.

The paper is organized as follows: In the section 1, we recall 
basic concepts and results in defoemation quantization, 
$\lambda$-formal analytic functions and the Gutt star product on 
$\Geh$. 
The main results of this paper are contained in section 2. 
We provide $\exp_*(i\xi\Phi_*(X))$, an integral expression of 
$\Phi_*$ and proof of the differentiability of $\Phi_*$. 
In the section 3, we show the structure theorem of $G$-invariant 
star products on a coadjoint orbit.

\section{Preliminarie}
\subsection{Star products}

Let $(M,\omega)$ be a symplectic manifold and $C^\infty(M)$ the set of 
smooth functions on $M$. 
The Poisson bracket on $C^\infty(M)$ associated to $\omega$ is denoted 
by $\{\ ,\ \}$. Let $C^\infty(M)[[\lambda]]$ be the space 
of power series in a formal parameter $\lambda$ with coefficients 
in $C^\infty(M)$.

A (differentiable) star product is an associative 
multiplication $*$ on $C^\infty(M)[[\lambda]]$ of the form
\begin{align*}
u*v=uv+\sum_{n=1}^\infty\left(\frac{\lambda}{2}\right)^nC_n(u,v),
\quad\text{for any }u,v\in C^\infty(M),
\end{align*}
where each $C_k$ is a bidifferential operator annihilating constants 
and $C_1(u,v)-C_1(v,u)=2\{u,v\}$. 
In the situation where a Lie group $G$ acts on $M$, a star product 
$*$ is said to be $G$-invariant if $g(u*v)=gu*gv$ 
holds for any $u,v\in C^\infty(M)[[\lambda]]$ and $g\in G$, 
where $gu(x)=u(g^{-1}x), x\in M$.
There exists a star product on any symplectic manifold 
\cite{WL:1983,OMY:1991,Fedosov:1994}, 
and the existence of $G$-invariant star products is equivalent to 
the existence of a $G$-invariant connection on $M$ 
\cite{Xu:1998,Fedosov:1996}. 
When $M$ is compact, $G$-invariant connections always exist and 
consequently there always exist $G$-invariant star product on $M$. 

Two star products $*_1$ and $*_2$ on $C^\infty(M)[[\lambda]]$ are said to be 
formally equivalent if there is a formal series,
\begin{align*}
T=Id+\sum_{n=1}^\infty \lambda^n T_n,
\end{align*}
of differential operators on $C^\infty(M)$ annihilating constants 
such that $u*_2v=T(T^{-1}u*_1T^{-1}v)$.
In this case, $T$ is called an equivalence between $*_1$ and $*_2$, 
and $*_2$ is denoted by $*_1^T$.
If $*_1$ and $*_2$ are equivalent  $G$-invariant star products and 
if the equivalence 
$T$ is $G$-invariant, then $*_1$ and $*_2$ are said to be formally 
$G$-equivalent and $T$ is called $G$-equivalence; 
see also \cite{BBG:1997,BBG:1998}.

\subsection{Formal analytic functions}
In this subsection, we make some simple, but useful remarks on the 
convergence of the power series valued in $\C[[\lambda]]$, 
that will be needed for calculus of functions in 
$C^\omega(\R^n)[[\lambda]]$.

\begin{defi}
A function $u=u_0+\lambda u_1+\cdots \in C^\infty(\R^n)[[\lambda]]$ 
is called formal analytic if each $u_i$ is analytic on $\R^n$.
We denote the set of formal analytic functions by 
$C^\omega(\R^n)[[\lambda]]$.
\end{defi}

Let $u$ and $v$ be formal analytic functions. 
We shall define the composition $u(v)$. 
If $u$ is a polynomial, there is no difficulty to define it, 
that is, it is given by substituting $v$ in $u$. 
For the general case, we define the composition by using power 
series. 
We begin with the following definition.

\begin{defi}
Let $a_J=\sum_{k=0}^\infty a_{J,k}\lambda^k\in \C[[\lambda]]$ 
be a multi-indexed sequence with respect to $J=(j_1,\ldots,j_n)$. 
The series $\sum_J a_J$ 
is said to converge formally absolutely if for any k, the series 
$\sum_J a_{J,k}$ converges absolutely.
\end{defi}
If a power series $\sum_J a_J y^J$ converges formally absolutely 
for some radius $\rho>0$ then this power series defines a formal 
analytic function on $|y|<\rho$.

Let $p^j(x)=\sum_{k=1}^\infty p^j_k(x)\lambda^k:
\R^m\to \lambda \R^n[[\lambda]]$, 
$j=1,\ldots n$, be a formal analytic map. 
A formal differential operator $p\partial$ is defined by 
$((p\partial)u)(y)=\sum p_k^j(x)(\partial_j u)(y)$ for $u:\R^n \to \R$. 
We define a formal operator $e^{ p\partial}$ for $u$ by
\begin{align}
(e^{ p\partial}u)(y) = 
u(y)+ \sum_{0<|J|}\frac{1}{|J|!}p^J(x)(\partial_J u)(y) .
\label{eq:ExtTaylor}
\end{align}
We should note that the right-hand side converges with respect to the 
filtration of $\lambda$ since $\deg(p)>0$. 
It is easy to show that $e^{ p\partial}u$ is an automorphism, that is, 
$e^{ p\partial}(u_1u_2)=(e^{ p\partial}u_1)( e^{ p\partial}u_2)$.

If $u$ is a polynomial on $\R^n$ then $u(y+p(x))$ is a function of 
$(x,y)$ that can be defined by substituting $y+p(x)$ in $u$, 
and we have $u(y+p(x))=(e^{ p\partial}u)(y)$. 
If $u$ is given by a power series, $u(y)=\sum_J a_J y^J$, 
one can see that the series $\sum_J a_J(y+p(x))^J$ is 
equal to $(e^{ p\partial}u)(y)$ as formal power series in $y$. 
Since $(e^{ p\partial}u)(y)$ converges formally absolutely on the 
same domain of $y\in \R^n$ where $\sum_J a_J y^J$ converges, 
we can define $u(y+p(x))=\sum_J a_J(y+p(x))^J$ as a formal analytic 
function. 
Therefore we can define $u(v)$ for any formal analytic map 
$v:\R^m\to \R^n[[\lambda]]$, and $u(v)=(e^{ (v-v_0)\partial}u)(v_0)$ 
holds, where $v=v_0+v_1\lambda+\cdots$.

\begin{defi}
Assume that $u:\R^n\to\R$ is an analytic map and let
$v:\R^m \to \R^n[[\lambda]]$ be a formal analytic map. 
Then we define a formal analytic map 
$u(v):\R^m\to\R[[\lambda]]$ by the following power series
\begin{align*}
u(v(x))=\sum_J a_J(v(x))^J,
\end{align*}
where $u=\sum_J a_J y^J$. 
\end{defi}

\remark \ The equation
\begin{align}
u(v(x))=(e^{ (v-v_0)(x)\partial}u)(v_0(x))
\label{eq:composition}
\end{align}
holds for any formal analytic map, and it gives the Taylor theorem 
for formal analytic functions.

\subsection{The Gutt star product}

Let $\Geh$ be a real Lie algebra and $\Geh^*$ its dual. 
The universal enveloping algebra (resp. symmetric algebra) of 
$\Geh$ is denoted by $\U(\Geh)$(resp. $\S(\Geh)$). 
We also denote the space of polynomials on $\Geh^*$  by $\pol(\Geh^*)$.
Let $\Geh[[\lambda]]$ be the set of formal power series in 
$\lambda $ with coefficients in $\Geh$. 
We define a Lie algebra structure $[\ ,\ ]_{\lambda }$ on 
$\Geh[[\lambda]]$ by 
$[\xi ,\eta ]_{\lambda }=\lambda \lbrack \xi ,\eta ]$ 
for any $\xi ,\eta \in \Geh$ and extend it by $\lambda$-linearity, 
where $[\ ,\ ]$ is the Lie bracket of $\Geh$. 
We denote this Lie algebra  by $\Geh_{\lambda }$. 
One can introduce a grading on $\Geh_\lambda$ by assigning 
to $\xi\in\Geh$, $\deg(\xi)=2$ and $\deg(\lambda)=2$, and 
$[\ ,\ ]_\lambda$ has the degree 0. 
This grading induces a grading on the universal enveloping 
algebra $\U(\Geh_\lambda)$ of $\Geh_\lambda$.

It is well known that the space of smooth functions on $\Geh^*$ 
admits a natural Poisson structure defined by the Kirillov-Poisson 
bracket $\Pi$.
For any smooth functions $u$ and $v$ on $\Geh^*$, $\Pi$ is given by
$\Pi(u,v)(\mu)=\langle [du(\mu),dv(\mu)],\mu \rangle,$
where $du(\mu)$ is an element of $\Geh$ considered as $1$-form on 
$\Geh^*$.

S.\ Gutt has defined a star product on $\Geh^*$\ \cite{Gutt:1983}. 
We shall call this product the Gutt star product, denoted by $*^G$. 
The Gutt star product can be directly obtained by transposing the 
algebraic structure of $\U(\Geh_\lambda)$ to 
$C^\infty(\Geh^*)[[\lambda]]$. 
This is achieved through the natural isomorphism between 
$\pol(\Geh^*)[[\lambda]]$ and $\S(\Geh_\lambda)$ and 
with the help of the symmetrization map 
$s:\S(\Geh_\lambda)\to\U(\Geh_\lambda)$. 
For polynomials $u$ and $v$, the Gutt star product is given by
\begin{align}
u *^G v=s^{-1}(s(u)\cdot s(v)),
\label{eq:defGutt}
\end{align}
where $\cdot$ is the product of $\U(\Geh_\lambda)$. 
Formula (\ref{eq:defGutt}) defines an associative differentiable 
deformation of the usual product on $\pol(\Geh^*)$. 
which admits a unique extension to $C^\infty(\Geh^*)[[\lambda]]$

As a direct consequence of Equation (\ref{eq:defGutt}) $*^G$ 
is a Weyl star product, that is, for any linear function 
$\xi$ on $\Geh^*$, $\xi^{*^G k}=\xi^k$ holds,
where $\xi^{*^G k}=\xi*^G\cdots *^G\xi$ ($k$ factors).
 Moreover, $*^G$ is $\Geh$-covariant,
\begin{align*}
\xi *^G \eta-\eta *^G \xi=
2\lambda \Pi(\xi,\eta)\quad \text{for }\xi,\eta\in \lin(\Geh^*),
\end{align*}
and $\Ad^*(G)$-invariant,
\begin{align*}
g(u*^G v)=(gu)*^G(gv)\quad 
\text{for } u,v\in C^\infty(\Geh^*)[[\lambda]],\ g\in G.
\end{align*}

There is a characterization of the Gutt star product given 
by the following.

\begin{prop}
[\protect\cite{Dito:1999}]
The Gutt star product is the unique $\Geh$-covariant Weyl star 
product on $(\Geh^*,\Pi)$. Any $\Geh$-covariant star product on 
$(\Geh^*,\Pi)$ is equivalent to the Gutt star product.
\end{prop}

Let $\xi=\sum_{k=0}^\infty \xi_k \lambda^k \in 
\lin(\Geh^*)[[\lambda]] \cong \Geh_\lambda$. 
Then a power series 
\begin{align}
e^\xi=\sum_{k=0}^\infty \frac{1}{k!}\xi^k \label{form:exp1}
\end{align}
can be defined in the sense of formal absolute convergence, and satisfies 
Equation (\ref{eq:composition}) . A simple computation implies 
that there are polynomials  $p_k(\xi_1,\xi_2,\cdots,\xi_k)$ such that 
\begin{align}
e^\xi=e^{\xi_0} \sum_{k=0}^\infty p_k(\xi_1,\xi_2,\cdots,\xi_k)\lambda^k.
\label{eq:starExpG}
\end{align}
Since $*^G$ is a Weyl star product, we also have 
\begin{align*}
e^\xi=\exp_{*^G}(\xi) \equiv \sum_{k=0}^\infty \frac{1}{k!}\xi^{*^G k}.
\end{align*}

For any $\xi,\eta\in\Geh_\lambda$, we denote by $\CH_\lambda(\xi,\eta)$ the
Campbell-Hausdorff series of a Lie algebra $\Geh_\lambda$.
We should note that $\CH_\lambda(\xi,\eta)$ is an element of $\Geh_\lambda$ 
since $[\ ,\ ]_\lambda$ has the degree $0$ and $\CH_\lambda$ converges with 
respect to the filtration of $\lambda$.

Since $*^G$ is $\Geh$-covariant, we have 
$\exp_{*^G}(\xi) *^G \exp_{*^G}(\eta)=\exp_{*^G}(\CH_\lambda(\xi,\eta))$, 
that is,
\begin{align}
e^{\xi} *^G e^{\eta}=
 e^{\CH_\lambda(\xi,\eta)} \quad\text{for }\xi,\eta\in\Geh_\lambda. 
\label{eq:expGroup}
\end{align}
Therefore, the set 
$G_\lambda \equiv \{e^\xi; \xi\in\Geh_\lambda \} \subset 
C^\omega(\Geh^*)[[\lambda]]$ is closed under the multiplication 
by $*^G$. 
It is also easy to show that $G_\lambda$ is a group.

\subsection{Oscillatory integral formula for star products}
For later use, we provide {\em an oscillatory integral expression} 
of the Gutt star product. 
We shall use the notations and the definitions given in the 
Appendix for oscillatory integrals; see also \cite{Kumano-go:1981}.

\begin{defi}
A function $u\in C^\infty(\R^n)$ has polynomial increase of degree 
$\tau>0$ if for any multi-index $I=(i_1,i_2,\cdots,i_n)$, there is a 
constant $C_I$ such that
\begin{align*}
|\partial_\zeta^I u(\zeta)|\leq C_I \langle \zeta \rangle^\tau.
\end{align*}
We can identify the set of these functions with $A_\tau^0$. 

\end{defi}
Let $\A^0=\bigcup_{\tau\geq 0}\A_\tau$. 
If we identify $\Geh^*$ to $\R^n$, $\A^0[[\lambda]]$ is a 
subalgebra of $(C^\infty(\Geh^*)[[\lambda]],*^G)$, which contains 
all polynomials.

\begin{defi}
Let $u\in \A^0$. 
The oscillatory integral expression of $u$ is given by 
the following formula,
\begin{align*}
u(\zeta)=\osint e^{i\alpha(\zeta-\beta)}u(\alpha)d\alpha d\beta,
\end{align*}
where the right-hand side means oscillatory integral.
\end{defi}
Since $*^G$ is differentiable, the $*^G$ operation commutes with 
integration. 
Therefore we have the oscillatory integral expression of 
the Gutt star product as follows: for any $u,v\in A^0$,
\begin{align*}
u *^G v(x)
&=\osint e^{-i(\alpha\beta+\alpha'\beta')}f(\alpha)g(\alpha')
 e^{i\beta x} *^G e^{i\beta' x}d\alpha d\alpha' d\beta d\beta'\\
&=\osint e^{-i(\alpha\beta+\alpha'\beta')}f(\alpha)g(\alpha')
 e^{\CH_\lambda(i\beta x,i\beta' x)}d\alpha d\alpha' d\beta d\beta'.
\end{align*}
We should remark that 
$e^{\CH_\lambda(i\beta x,i\beta' x)}\in 
\A^0[[\lambda]]\times \A^0[[\lambda]]$ 
because of Equation (\ref{eq:starExpG}). 
Hence the above computation makes sense and $\A^0[[\lambda]]$ 
is a subalgebra of $(C^\infty(\Geh^*)[[\lambda]],*^G)$.


\section{Differentiability of quantum moment map}


This section is devoted to the study of quantum moment maps 
a main subject of this paper. 
The definition of quantum moment maps, which we adopt here, 
is given in \cite{Xu:1998}; see Definition \ref{def:QMM}. 
This definition is a natural analogue of the definition of 
classical moment maps in Hamiltonian systems. 

However quantum moment maps are different from their classical 
counterparts in a significant way, by the locality features.
In the classical case, to give a ring morphism of $C^\infty(\Geh^*)$ 
into $C^\infty(M)$ is equivalent to give a differential map of 
$M$ into $\Geh^*$, that is a consequence from the locality of the 
ring of functions and its ring morphisms. 
So this implies that any ring morphism of $\pol(\Geh^*)$ into 
$C^\infty(M)$ has a natural extension to $C^\infty(\Geh^*)$. 
But the problem is not clear for the quantum case. 
There is no guarantee that a homomorphism of star algebras is local 
or differentiable. 

We shall here show that any quantum moment map is differentiable.

\subsection{Definition of quantum moment maps}

Let $(M,\omega)$ be a symplectic $G$-space and $*$ a $G$-invariant 
star product. 
We denote the star commutator by $[a,b]_*=a*b-b*a$. 
\begin{defi}
[\cite{Xu:1998}]
\label{def:QMM}
A quantum moment map is a homomorphism of associative algebras
\begin{align}
\Phi _* :\mathfrak{U}({\mathfrak{g}}_{\lambda })\rightarrow 
C^\infty(M)[[\lambda ]],
\label{eq:qmm1}
\end{align}
which satisfies 
\begin{align}
[ \Phi _*(\xi),u]_* =\lambda \xi u,  \label{eq:qmm2}
\end{align}
where the right-hand side of (\ref{eq:qmm2}) is the infinitesimal action
of $\xi$ $\in \Geh$ on $C^\infty(M)[[\lambda ]]$. 
\end{defi}
It is easy to see that the condition (\ref{eq:qmm1}) is equivalent to
\begin{align}  \label{eq:qmm3}
\Phi _*([\xi,\eta]_{\lambda })=[\Phi _*(\xi),\Phi _*(\eta)]_{ * }\quad \text{ for any }
\xi,\eta\in \Geh.
\end{align}

On the existence and the uniqueness of quantum moment maps, some simple
criteria are known.

\begin{thm}
[\protect\cite{Xu:1998}] 
\label{thm:extQMM}
Let $\text{H}_{dR}^*(M)$ be the de Rham cohomology group and 
$\text{H}^*(\Geh,\ \R)$ be Lie algebra cohomology group with
coefficients in $\R$. 
There exists a quantum moment map if 
$\text{H}_{dR}^{1}(M)=0$ and $\text{H}^{2}(\Geh,\ \R)=0$.
\end{thm}

\begin{thm}
[\protect\cite{Xu:1998}]
\label{thm:uniQMM}
 The set of quantum moment maps of a $G$-invariant star
product is parametrized by $\text{H}^1(\Geh,\R)$.
\end{thm}

The following proposition says that a quantum moment map is a natural 
analogue of the classical. 
\begin{prop}
[\cite{Xu:1998}]
Let $\Phi_*:\pol(\Geh^*[[\lambda ]])\rightarrow C^\infty(M)[[\lambda ]]$ 
be a quantum moment map. Then $M$ is a Hamiltonian $G$-space. 
Moreover $\Phi_*$ satisfies
\begin{equation*}
\Phi _*(u)=\Phi_0 (u)+O(\lambda ),\text{ for any }u \in 
\pol(\Geh^*),
\end{equation*}
where $\Phi_0 :\pol(\Geh^*)\rightarrow C^\infty(M)$
denotes the corresponding classical moment map.
\end{prop}

An important property of $\Phi_*$ is its covariance 
under $G$-equivalence.
\begin{prop}
Let $*$ be a $G$-invariant star product and $\Phi_*$ a quantum 
moment map of $*$. 
If $*'$ is a $G$-invariant star product which is $G$-equivalent to 
$*$, then $T\Phi_*$ is a quantum moment map of $*'$, where $T$ is 
a $G$-equivalence between $*$ and $*'$.
\end{prop}
\begin{proof}\label{prop:qmmTrans}
It is enough to show that $[T\Phi_*(X),f]_{*'}=\lambda Xf$ since 
$T\Phi_*$ is an algebra homomorphism from 
$(\pol(\Geh^*)[[\lambda]],*^G)$ to $(C^\infty(M)[[\lambda]],*')$.
\begin{align*}
[T\Phi_*(\xi),f]_{*'}=T[\Phi_*(\xi),T^{-1}f]_*
=T(\lambda \xi T^{-1}f)=\lambda \xi f.
\end{align*}
\end{proof}
Since quantum moment maps are parametrized by $H^1(\Geh,\R)$ we have,
\begin{cor}
Assume $H^1(\Geh,\R)=\{0\}$. 
Let $*,*'$ be $G$-invariant star product and $\Phi_*,\Phi_{*'}$ 
quantum moment maps of $*,*'$ respectively. 
If $*'$ is $G$-equivalent to $*$ then $T\Phi_*=\Phi_{*'}$.
\end{cor}

\subsection{Exponential function of a quantum moment map}

We shall define here a function $\exp_*(\Phi_*(X))$ in 
$C^\infty(M)[[\lambda]]$ for $X\in \Geh_\lambda$. 
This function ``generates''  $\Phi_*(\pol(\Geh^*))$, and we will 
use it to obtain another expression for $\Phi_*$. 
An important property of $\exp_*(\Phi_*(X))$ is that 
it is a product of $e^X$ and 
a polynomial of $X$.
This property is essential for the differentiability of $\Phi_*$.

Assume that $*$ is a $G$-invariant star product of Fedosov type. 
Recall that $Q$ and $\sigma$ denote the Fedosov quantization 
procedure and the projection of $W_D$ onto $C^\infty(M)[[\lambda]]$ 
respectively.

\begin{lem}\label{lem:starExp}
Let $\xi=\xi_0+\xi_1\lambda+\cdots \in\Geh_\lambda$. 
The series
\begin{align}
\sum_{k=0}^\infty \frac{1}{k!}Q(\Phi_*(\xi))^{\circ k}
\label{form:exp2}
\end{align}
converges $\lambda$-formally absolutely and uniformly on any compact 
subset of $M$, and defines an element of $\Gamma W_D$. 
Moreover, (\ref{form:exp2}) has the following expression
\begin{align}
\sum_{k=0}^\infty \frac{1}{k!}Q(\Phi_*(\xi))^{\circ k}=
e^{\Phi_0(\xi_0)}\sum_{I,j} 
 p_{I,j}(\Phi_*(\xi),\partial\Phi_*(\xi),\cdots)y^I\lambda^j ,
\end{align}
where $p_{I,j}$ are polynomials in 
$\{\Phi_*(\xi),\partial\Phi_*(\xi),\cdots\}$.
\end{lem}

\begin{proof}
If we decompose $Q(\Phi_*(\xi))=\Phi_0(\xi_0)+R(\xi)$, 
where $\Phi_0$ is the classical moment map, then $\deg\ R(\xi)\geq 1$ 
and $[\Phi_0(\xi_0),R(\xi)]_\circ =0$, 
since $\Phi_0(\xi_0)$ is a central element of $\Gamma W$. 
Therefore we have the following (formal) equation: 
\begin{align}
\sum_{k=0}^\infty \frac{1}{k!}Q(\Phi_*(\xi))^{\circ k}=
\sum_{k=0}^\infty \frac{1}{k!}\Phi_0(\xi_0)^k 
\sum_{k=0}^\infty \frac{1}{k!}R(\xi)^{\circ k}.
\label{eq:exp2}
\end{align}
The second factor of the r.h.s in (\ref{eq:exp2}) converges with respect 
to the filtration of $\Gamma W$ since $\deg\ R(\xi)\geq 1$.
So, it is easy to see that the r.h.s of (\ref{eq:exp2}) converges 
absolutely and uniformly on any compact subset of $M$.
Applying the Weyl derivation $D$ on (\ref{form:exp2}) 
term by term, we see that (\ref{form:exp2}) is a flat section.

To show the last statement, we express $R(\xi)$ as follows
\begin{align}
R(\xi)=\sum_{|I|+j \geq 1} 
r_{I,j}(\Phi_*(\xi),\partial\Phi_*(\xi),\cdots) y^I \lambda^j.\label{eq:rFQ}
\end{align}
Each $r_{I,j}$ is a polynomial in $\{\Phi_*(\xi),\partial\Phi_*(\xi),\cdots\}$ 
since it is obtained by Fedosov quantization procedure. 
Therefore each coefficient of $y^I \lambda^j $ in the series
\begin{align}
\sum_{k=0}^\infty \frac{1}{k!}R(\xi)^{\circ k} \label{eq:rExp2}
\end{align}
is also a polynomial.
\end{proof}

\begin{defi}
For any $\xi\in\Geh_\lambda$, the function $\exp_*(\Phi_*(\xi))$ in 
$C^\infty(M)[[\lambda]]$ is defined by
\begin{align*}
\exp_*(\Phi_*(\xi))
=\sigma(\sum_{k=0}^\infty \frac{1}{k!}Q(\Phi_*(\xi))^{\circ k}).
\end{align*}
\end{defi}

In the proof of Lemma \ref{lem:starExp}, assigning $\xi=\alpha^l X_l$, 
where $\{\alpha^l\}\in\C^n[[\lambda]]$ and $\{X_l\}$ is a basis of 
$\Geh$, we have the following.
\begin{cor}\label{cor:starExp1}
$\exp_*(\Phi_*(\alpha^l X_l))$ is a product of 
$e^{\alpha^l_0 \Phi_0(X_l)}$ and a polynomial of $\alpha^l$ taking 
values in $C^\infty(M)[[\lambda]]$.
\end{cor}

\begin{proof}
Since a quantum moment map is linear with respect to 
$\xi=\alpha^l X_l\in\Geh[[\lambda]]$, then each $r_{I,j}$ of 
(\ref{eq:rFQ}) is also linear with respect to $\xi$ and 
(\ref{eq:rExp2}) is a polynomial in $\alpha^l$.
\end{proof}

\begin{lem}\label{lem:starExpForm1}
Assume $\xi,\eta\in \Geh_\lambda$. Then we have 
\begin{align}
\exp_*(\Phi_*(\xi))*\exp_*(\Phi_*(\eta))=\exp_*(\Phi_*(\CH_\lambda(\xi,\eta))).
\end{align}
\end{lem}

\begin{proof}
By the definition of $\exp_*(\Phi_*(\xi))$,
\begin{align*}
Q(\exp_*(\Phi_*(\CH_\lambda(\xi,\eta)))
&=\sum_{k=0}^\infty
	\frac{1}{k!}Q(\Phi_*(\CH_\lambda(\xi,\eta)))^{\circ k}\\
&=\sum_{k=0}^\infty
	\frac{1}{k!}(\CH_\circ(Q(\Phi_*(\xi)),Q(\Phi_*(\eta))))^{\circ k},
\end{align*}
where 
$\CH_\circ$ denotes the Campbell-Hausdorff series with respect to 
the Weyl product $\circ$ of Weyl bundle $\Gamma W$. 
Since the equation 
\begin{align*}
\sum_{k=0}^\infty 
	\frac{1}{k!}(\CH_\circ(Q(\Phi_*(\xi)),Q(\Phi_*(\eta))))^{\circ k} =
\sum_{k=0}^\infty
	\frac{1}{k!}(Q(\Phi_*(\xi)))^{\circ k} \circ 
\sum_{k=0}^\infty
	\frac{1}{k!}(Q(\Phi_*(\eta)))^{\circ k}.
\end{align*}
holds in $\Gamma W$, we have the lemma. 
\end{proof}

With each multi-index $J=(j_1,j_2,\cdots,j_n)$ we associate 
a differential operator 
\begin{align*}
D^J_\alpha=\left(-i\frac{\partial}{\partial \alpha^1}\right)^{j_1}
\cdots
\left(-i\frac{\partial}{\partial \alpha^n}\right)^{j_n}.
\end{align*}

\begin{lem}
Assume $\{\alpha^l\}\in\R^n$. Then we have 
\begin{align}
(D^J \exp_*(\Phi_*(i \alpha^l X_l)))|_{\alpha=0}=\Phi_*(X^J).
\label{eq:genRule}
\end{align}
\end{lem}
\begin{proof}
Let $\tilde{X}_l=Q(\Phi_*(X_l))$.
\begin{align*}
D^J(Q(\exp_*(\Phi_*(i \alpha^l X_l))))|_{\alpha=0}
&=D^J(\sum_{k=0}^\infty \frac{1}{k!}(i\alpha^l \tilde{X}_l)^{\circ k})|_{\alpha=0}\\
&=D^J(\frac{1}{|J|!}(i\alpha^l \tilde{X}_l)^{|J|})|_{\alpha=0}\\
&=Q(\Phi_*(X^J)).
\end{align*}
\end{proof}

\subsection{Oscillatory integral expression for $\Phi_*$}

We shall provide another expression for a quantum moment map $\Phi_*$ 
by using $\exp_*(\Phi_*(X))$ and an oscillatory integral. 
This expression gives us a clear understanding of quantum moment maps 
and enables to show the differentiability of $\Phi_*$.

\begin{defi}
\label{def:extQmm}
Let $\{X_l\}$ be a basis of $\Geh$ and $\{X^l\}$ its dual basis. 
We define the map $\overline{\Phi}_*$ from 
$\A^0$ into $C^\infty(M)[[\lambda]]$ as follows. 
\begin{align}
\overline{\Phi}_*(u)=
\osint u(\mu X) e^{-i \nu \mu}\exp_*(\Phi_*(i \nu X)) 
d\mu d\nu, \quad u\in\A^0,
\end{align}
where $\mu X=\mu_l X^l$ and $\nu X=\nu^l X_l$.
\end{defi}
This definition makes sense since $u(\mu X)\exp_*(\Phi_*(i\nu X))\in \A$.
It is easy to see that the above definition does not depend on 
a choice of a basis $\{X_l\}$.
\begin{lem}
$\overline{\Phi}_*$ coincides with $\Phi_*$ on polynomials.
\end{lem}
\begin{proof}
Let $X^J$ be a monomial on $\Geh^*$. Then 
\begin{align*}
\overline{\Phi}_*(X^J)
&= \osint \beta^J e^{-i \alpha \beta}\exp_*(\Phi_*(i \alpha X)) d\alpha d\beta \\
&= \osint e^{-i \alpha \beta} D_\alpha^J\exp_*(\Phi_*(i \alpha X)) d\alpha d\beta \\
&= D_\alpha^J\exp_*(\Phi_*(i \alpha X))|_{\alpha=0}=\Phi_*(X^J),
\end{align*}
where we have applied Equation (\ref{eq:genRule}) in the last line.
\end{proof}

So we shall also use the notation $\Phi_*$ for $\overline{\Phi}_*$.

The following proposition says that $\exp_*(\Phi_*(X))$ can be 
considered as the image of $e^X$ under a quantum moment map.
\begin{prop}\label{prop:QMMExtForm1}
Assume that $p^k=p^k_j\lambda^j\in\C[[\lambda]]$ satisfies 
$p^k_0\in i\R$.
Then $e^{pX}=e^{p^kX_k}\in \A$ and $\Phi_*(e^{pX})=\exp_*(\Phi_*(pX))$.
\end{prop}

\begin{proof}
Let $pX=ia X+r X$, where $a^k\in \R$ and 
$r\in \lambda\C[[\lambda]]$. Then $e^{pX}\in \A$ by Equation 
(\ref{eq:starExpG}). 
By Definition of $\Phi_*$, 
\begin{align*}
\Phi_*(e^{px})
&=\osint e^{p\mu} e^{-i\mu\nu} \exp_*(\Phi_*(i\nu X)) d\mu d\nu \\
&=\osint e^{i a \mu} e^{r \mu} e^{-i\mu\nu} \exp_*(\Phi_*(i\nu X)) d\mu d\nu \\
&=\osint e^{i (a-\nu)\mu} (e^{r D_\nu} \exp_*(\Phi_*(i\nu X))) d\mu d\nu \\
&=\osint e^{i (a-\nu)\mu} \exp_*(\Phi_*(i(\nu-ir) X)) d\mu d\nu \\
&=\exp_*(\Phi_*(i(a-ir) X))=\exp_*(\Phi_*(pX)),
\end{align*}
where we have applied Equation (\ref{eq:composition}).
\end{proof}

As a corollary of Proposition \ref{prop:QMMExtForm1} and Lemma 
\ref{lem:starExpForm1}, we have 
$\Phi_*(e^{i \xi})*\Phi_*(e^{i \eta})=\Phi_*(e^{i \xi} *^G e^{i \eta})$.

\begin{thm}
\label{thm:qmmDiff}
A quantum moment map $\Phi_*$ is differentiable. 
Moreover, if $*$ is of Fedosov type, then there are functions 
$S_{I,j}\in C^\infty(M),\ I=(i_1,\ldots,i_{n})$, 
$j=0,1,\ldots$ such that 
\begin{align}
\Phi_*(u)=\sum_{j=0}^\infty \lambda^j \sum_{\ 0 \leq |I| \leq 2j}
 S_{I,j}\Phi_0(D_\mu^I u), 
\quad\text{for any }u(\mu)\in C^\infty(\Geh^*). 
\end{align}
\end{thm}

\begin{proof}
First we assume that $*$ is of Fedosov type. 
By Corollary \ref{cor:starExp1}, we have an expression for 
$\exp_*(\Phi_*(i\alpha^l X_l))$ given by 
\begin{align*}
e^{i\alpha^l \Phi_0(X_l)}\sum_{j=0}^\infty \lambda^j \sum_{0 \leq |I| \leq 2j}
S_{I,j}\alpha^I ,
\end{align*}
where $S_{I,j}\in C^\infty(M)$ depends on $\Phi_*(X_i)$ and on the 
Weyl connection $D$ of $*$. 
By the definition of $\Phi_*$, we have 
\begin{align*}
\Phi_*(u)
&=\osint u(\mu)e^{-i\mu_l\nu^l}\exp_*({\Phi_*(i\nu^l X_l)}) d\mu d\nu \\
&=\sum_{j=0}^\infty \lambda^j \sum_{0 \leq |I| \leq 2j} 
\osint u(\mu)e^{-i\mu_l\nu^l}e^{i\nu^l\Phi_0(X_l)} S_{I,j}\nu^I d\mu d\nu \\
&=\sum_{j=0}^\infty \lambda^j \sum_{0 \leq |I| \leq 2j} 
\osint u(\mu)e^{-i\nu^l(\mu_l-\Phi_0(X_l))} S_{I,j}\nu^I d\mu d\nu \\
&=\sum_{j=0}^\infty \lambda^j \sum_{0 \leq |I| \leq 2j} 
\osint S_{I,j}(D_\mu^I u)(\mu)e^{-i\nu^l(\mu_l-\Phi_0(X_l))} d\mu d\nu \\
&=\sum_{j=0}^\infty \lambda^j \sum_{0 \leq |I| \leq 2j} S_{I,j}\Phi_0(D_\mu^I u).
\end{align*}
For a general $G$-invariant star product $*'$, we have a $G$-invariant 
star product $*$ of Fedosov type which is $G$-equivalent to $*'$ 
by Theorem \ref{thm:GStarClass}. 
If $T$ denotes a $G$-equivalence between $*$ and $*'$, then 
any quantum moment map $\Phi_{*'}$ of $*'$ has the form $T\Phi_*$. 
Therefore $\Phi_{*'}$ is differentiable.
\end{proof}

As a consequence of Theorem \ref{thm:qmmDiff}, $\Phi_*$ admits a unique 
extension to $C^\infty(\Geh^*)[[\lambda]]$. 
Since $\Phi_*$ is an algebra homomorphism on polynomials on $\Geh^*$, 
the differentiability of $\Phi_*$ implies $\Phi_*$ 
is an algebra homomorphism on $C^\infty(\Geh^*)[[\lambda]]$.

It is not difficult to compute $S_{I,j}$ for lower degrees 
in ${I,j}$. 
For instance, we have 
$S_{0,0}=1,S_{l,1}=\Phi_1(X_l),S_{lm,1}
=\{\Phi_0(X_l),\Phi_0(X_m)\}$, and so on. 
Therefore we have $\Phi_*(u)=\Phi_0(u)+o(\lambda)$ for any 
$u\in C^\infty(\Geh^*)$.

\subsection{Properties of $\Phi_*$ }

\begin{prop}
\label{prop:gcov}
A quantum moment map is a $\Geh$-equivariant map from 
$C^\infty(\Geh^*)[[\lambda]]$ to $C^\infty(M)[[\lambda]]$. 
Therefore, $\Phi_*$ is $G$-equivariant if $G$ is connected.
\end{prop}
\begin{proof}
This proposition is a direct consequence from the definition of quantum 
moment maps and $G$-invariance of $*^G$ and $*$; For any 
$u\in C^\infty(\Geh^*)$ and $\xi\in\Geh$,
\begin{align*}
\Phi_*(\lambda \xi u)=\Phi_*([\xi,u]_{*^G})=[\Phi_*(\xi),\Phi_*(u)]_*=
\lambda \xi\Phi_*(u).
\end{align*}
\end{proof}

\begin{prop}\label{prop:qmmGInv}
If $f\in C^\infty(M)[[\lambda]]$ commutes with any 
$\Phi_*(u),u\in C^\infty(\Geh^*)$ then $f$ is a $G$-invariant function.
\end{prop}
\begin{proof}
It is easy.
\end{proof}

\begin{prop}\label{prop:qmmSur}
A quantum moment map $\Phi_*$ is surjective if and only if its classical part 
$\Phi_0$ is surjective.
\end{prop}
\begin{proof}
Assume that $\Phi_0$ is surjective.
For $u=\sum u_i\lambda^i\in C^\infty(\Geh^*)[[\lambda]]$ and 
$\varphi\in C^\infty(M)$ the equation $\Phi_*(u)=\varphi$ is 
equivalent to
\begin{align}
\Phi_0(u_0)&=\varphi, \label{eq:qmmSur1} \\
\Phi_0(u_k)&=-\sum_{j=1}^{k}\Phi_j(u_{k-j})\quad\text{for any }k>0.
\label{eq:qmmSur2}
\end{align}
One can solve the above system of equations by induction, since $\Phi_0$ 
is surjective. The converse is trivial.
\end{proof}

\begin{lem}\label{lem:qmmSolLocal}
Let $\varphi$ be a smooth function on $M$. 
Assume there exists $u\in C^\infty(\Geh^*)[[\lambda]]$, 
that is a solution of the equation $\Phi_*(u)=\varphi$. 
Then $u$ depends locally on $\varphi$. 
More precisely, the dependence of $u$ at $J(q)$ on $\varphi$ 
is described by differentials of $\varphi$ at $q \in M$, 
where $J:M\to\Geh^*$ is the dual form of $\Phi_*$, that is, 
$(\Phi_0(u))(q)=u(J(q))$.
\end{lem}
\begin{proof}
Equation (\ref{eq:qmmSur1}) says that $u_0(J(q))$ depends on 
$\varphi(q)$. Since $\Phi_*$ is differentiable, the right-hand side of 
Equation (\ref{eq:qmmSur2}) also depends on differentials of $\varphi$  
if $u_0,\ldots,u_{k-1}$ depend on differentials of $\varphi$.
\end{proof}


\subsection{Invariants for $G$-invariant star products on transitive spaces}


In this subsection we review the results of 
\cite{Hamachi:1999,Hamachi:2002}, where we define invariants 
for $G$-invariant star products on $G$-transitive symplectic manifolds.

Let $M$ be a $G$-transitive symplectic manifold and $*$ a 
$G$-invariant star product on $M$. 
We assume that there is a unique quantum moment map $\Phi_*$ for $*$. 
Let $\Z$ be the center of $C^\infty(\Geh^*)$, that is, 
the set of functions that commute with any smooth function on 
$\Geh^*$ with respect to the Gutt star product. 
One can show that $\Z$ is equal to the set of $G$-invariant functions 
on $\Geh^*$. For any $l\in\Z$ we have 
$[\Phi_*(l),\Phi_*(C^\infty(\Geh^*))]_*=\Phi_*([l,C^\infty(\Geh^*)]_{*^G})=0$, 
so that Proposition \ref{prop:qmmGInv} implies $\Phi_*(l)$ is a 
$G$-invariant function on $M$. 
Since $M$ is transitive, $\Phi_*(u)$ is constant. 
Consequently, we make the following definition.

\begin{defi}
Let $M$ be a $G$-transitive symplectic manifold and $*$ a $G$-invariant 
star product which admits a unique quantum moment map $\Phi_*$. 
A map $c_*$ from $\Z$ to $\C[[\lambda]]$ is defined by 
$c_*(l)\equiv \Phi_*(l)$ for any $l\in \Z$.
\end{defi}

The following simple proposition is important.
\begin{prop}
\label{prop:QMMKerC}
If $\ker \Phi_*=\ker \Phi_{*'}$ then $c_*=c_{*'}$.
\end{prop}
\begin{proof}
Any function on $C^\infty(\Geh^*)[[\lambda]]$ of the form $l-c_*(l)$ for 
$l\in\Z$ is an element of $\ker \Phi_*$. Therefore if 
$\ker \Phi_*=\ker \Phi_{*'}$, we have $\Phi_{*'}(l-c_*(l))=0$, that is, 
$c_{*'}(l)=c_*(l)$.
\end{proof}

\begin{cor} 
If $*'$ is $G$-equivalent to $*$ then $c_*=c_{*'}$.
\end{cor}

This Corollary means that $c_*$ depends only on the class of 
$G$-equivalent star products. 
We have computed $c_*$ for few a examples in \cite{Hamachi:2002}.

There is a natural question to ask:
``Does $c_*$ parameterize the class of $G$-invariant star products 
on a $G$-transitive space? '' 
In the next section, we give a complete answer of this question 
when $M$ is a coadjoint orbit of a compact semisimple Lie group.


\section{Star products on regular coadjoint orbits of compact semisimple Lie groups}


Let $G$ be a real compact semisimple Lie group and 
$\O \subset \Geh^*$ a regular coadjoint orbit of $G$. 
$\O$ has a natural symplectic structure that is induced from 
the Kirillov-Poisson structure $\Pi$. 

Recall that there is a $G$-invariant star product on $\O$ since 
$\O$ is compact. 
Since $G$ is semisimple, 
Theorems \ref{thm:extQMM} and \ref{thm:uniQMM} imply that 
for each $G$-invariant star product $*$ on $\O$, there is a 
unique quantum moment map of $*$.

We study here $G$-invariant star products on $\O$, and 
our goal is to present a structure theory for them.

\subsection{Structure theory}
Let $*$ be a $G$-invariant star product on $\O$ and 
$\Phi_*=\Phi_0+\Phi_1\lambda^1+\cdots$ a quantum moment map of $*$.
The classical moment map $\Phi_0$ is simply given by the 
pull-back of the embedding map of $\O$ in $\Geh^*$, 
that is, $\Phi_0$ is surjective. 
Therefore, Proposition \ref{prop:qmmSur} implies that the 
quantum moment map $\Phi_*$ is also surjective. 
The following proposition is a direct consequence. 
\begin{prop}\label{prop:homThm1}
We have a $G$-equivariant isomorphism 
\begin{align}
C^\infty(\Geh^*)[[\lambda]]/\ker \Phi_* \cong C^\infty(\O)[[\lambda]].
\end{align}
\end{prop}

Let $*$ and $*'$ be $G$-invariant star products. 
If we assume $\ker \Phi_*=\ker \Phi_{*'}$, 
Proposition \ref{prop:homThm1} defines a morphism 
$S:(C^\infty(\O)[[\lambda]],*) \to (C^\infty(\O)[[\lambda]],*')$. 
Let $\varphi$ be a smooth function on $\O$. 
There exists $u\in C^\infty(\Geh^*)[[\lambda]]$ such that 
$\Phi_*(u)=\varphi$, and we define $S(\varphi) \equiv \Phi_{*'}(u)$. 
Lemma \ref{lem:qmmSolLocal} and the differentiability of $\Phi_{*'}$ 
implies that the morphism $S$ is differentiable. 
Therefore we have the following lemma.
\begin{lem}
$\ker \Phi_*=\ker \Phi_{*'}$ if and only if $*$ is $G$-equivalent to $*'$.
\end{lem} 

As we have seen before, $\ker \Phi_*=\ker \Phi_{*'}$ implies $c_*=c_{*'}$. 
The following Proposition shows that if there are ``good coordinates'' 
on $\Geh^*$, the converse also holds.

\begin{prop}\label{prop:kerQMM}
Assume there are functionally independent $G$-invariant functions 
$p_i:\Geh^*\to \R,1\leq i \leq r$ such that $\O$ is given  by the 
level set $\{\xi \in \Geh^*:p_i(\xi)=c_i\}$ for some regular value 
$\{c_i\}$ of $\{p_i\}$. 
Then $\ker \Phi_*$ is equal to the ideal of 
$(C^\infty(\Geh^*)[[\lambda]],*^G)$ generated by 
$\{p_i-c_*(p_i);1\leq i \leq r\}$.
\end{prop}

\begin{proof}
Let $f=f_0+f_1\lambda+\cdots \in \ker \Phi_*$. 
Then it is easy to see that $f_0\in \ker \Phi_0$, that is, $f_0$ 
is null on $\O$. 
So there are functions $g_i\in C^\infty(\Geh^*)$ such that
\begin{align}
f_0=\sum_{i=1}^r g_i(p_i-c_i).
\end{align}
Setting 
\begin{align}
f^{(0)}=\sum_{i=1}^r g_i*^G(p_i-c_*(p_i)),\label{eq:genIdeal}
\end{align}
 we have $f^{(0)}\in \ker \Phi_*$ and 
$f^{(0)}_0=f_0$. Applying the same argument on $(f-f^{(0)})/\lambda$ 
inductively we find a sequence of functions $f^{(k)}$ satisfying 
\begin{align}
f=\sum_{k=0}^\infty f^{(k)}\lambda^k.
\end{align}
Since each $f^{(k)}$ has the form (\ref{eq:genIdeal}), 
it completes the proof.
\end{proof}

Let $I\subset \pol(\Geh^*)$ be the set of polynomials on $\Geh^*$ 
invariant under $G$. 
By Chevalley's theorem one has $I=\C[p_1,\ldots,p_r]$, where 
$p_1,\ldots,p_r$ are algebraically independent homogeneous 
polynomials and $r$ is the rank of $\Geh$ \cite{Va:1984}. 
One can also see that any regular coadjoint orbit $\O$ is given by 
the level set $\{\xi\in\Geh^*; p_1(\xi)=c_1,\ldots,p_r(\xi)=c_r\}$ 
for some regular value $\{c_j\}$ \cite{Kostant:1963}. 
Therefore, $\{p_j\}$ satisfies the condition of 
Proposition \ref{prop:kerQMM}. 
So we have the inverse of Proposition \ref{prop:QMMKerC}. 

\begin{prop}
\label{prop:kerQMMStr}
Let $*,*'$ be $G$-invariant star products and $\Phi_*,\Phi_{*'}$ be 
quantum moment maps of $*,*'$ respectively. 
Then $c_*=c_{*'}$ implies $\ker \Phi_*=\ker \Phi_{*'}$. 
Moreover, if we take $\{p_j\}$ as the algebraically independent 
homogeneous polynomials that are obtained by Chevalley's theorem, 
then we have 
\begin{align*}
\ker \Phi_*=\langle p_j-\Phi_*(p_j) \rangle,
\end{align*}
where $\langle p_j-\Phi_*(p_j) \rangle$ denotes the ideal of 
$C^\infty(\Geh^*)[[\lambda]]$ generated by $p_j-\Phi_*(p_j)$.
\end{prop}

And we have also the following structure theorem.

\begin{thm}
For any $G$-invariant star product $*$ on $\O$, there are constants 
$c_{*,j}\in \C[[\lambda]],j=1,2,\ldots,r$ such that 
\begin{align*}
(C^\infty(\O)[[\lambda]],*) \cong 
C^\infty(\Geh^*)[[\lambda]]/\langle p_i-c_{*,j} \rangle.
\end{align*}
Moreover, this isomorphism is $G$-equivariant. 
\end{thm}
\begin{proof}
Let $\Phi_*$ be a quantum moment map of $*$. 
If we set $c_{*,j}=\Phi_*(p_j)$, the theorem is a direct consequence 
of Propositions \ref{prop:homThm1} and \ref{prop:kerQMMStr}.
\end{proof}

\begin{cor}
There is a one-to-one correspondence between the classes of 
$G$-invariant star products on $\O$ and $c_*$. 
\end{cor}

\subsection{Example}
Coadjoint orbits of $\SO(3)$.

Let $\O$ be a regular coadjoint orbit of $\SO(3)$. 
It is well known that $\O$ is a two-dimensional sphere in 
$\so(3)^*$ given by the Casimir polynomial $p(x,y,z)=x^2+y^2+z^2$ 
on $so(3)^*=\R^3$,
that is, there is a real number $r>0$ such that
\begin{align*}
\O=\{(x,y,z)\in\Geh^*; p(x,y,z)=r^2\}.
\end{align*}

The class of $\SO(3)$-invariant star products is parametrized by 
the second equivariant de Rham cohomology $H_{dR}^2(\O,\R)^{\SO(3)}$. 
Let $*$ be a $\SO(3)$-invariant star product on $\O$. 
Then $p$ satisfies the condition of Proposition \ref{prop:kerQMM}, 
so that $\ker\Phi_*$ is described by $c_*(p)$  and we obtain
\begin{align}
(C^\infty(\O)[[\lambda]],*) \cong 
(C^\infty(\Geh^*)[[\lambda]],*^G)/\langle p-c_*(p) \rangle.
\label{eq:repEx1}
\end{align}
Hence, $G$-invariant star products on $\O$ have the form of 
the right-hand side of (\ref{eq:repEx1}) and are parametrized by 
$c_*(p)$. 
This gives another classification of $G$-invariant star products on 
$\O$.


\setcounter{section}{0}
\renewcommand{\thesection}{\Alph{section}}

\section*{Appendices}

\section{Oscillatory integral}


We provide here a brief review on {\em Oscillatory integrals} 
in order to fix some definitions and notations. 
The following is based on \cite{Kumano-go:1981}, with little 
modifications adapted to our problem.

\begin{defi}
A function $a(\xi,x)\in C^\infty(\R^n_\xi\times\R^n_x)$ is said to 
belong to the $\A_{\tau}^m$-class, $-\infty<m<\infty,\ 0\leq \tau$ 
if for any multi-indices $I$ and $J$, 
there exists a constant $C_{I,J}$ such that
\begin{align*}
|\partial_\xi^I \partial_x^J a(\xi,x)|
\leq C_{I,J}\langle \xi \rangle^m \langle x \rangle^\tau,
\end{align*}
where $\langle\xi\rangle=\sqrt{1+|\xi|^2}$. 
Set 
\begin{align*}
\A=\bigcup_{-\infty<m<\infty}\bigcup_{0\leq \tau} \A_\tau^m.
\end{align*}
\end{defi}
For $a(\xi,x)\in \A_\tau^m$, we define a family of seminorms 
$|a|_l,\ l=0,1,\cdots,$ by
\begin{align*}
|a|_l=\max_{|I+J|\leq l}\sup_{(\xi,x)}\{|\partial_\xi^I \partial_x^J a(\xi,x)| \langle \xi \rangle^{-m} \langle x \rangle^{-\tau} \}.
\end{align*}
Then $\A_\tau^m$ becomes a Fr\'echet space. A subset $B$ of $\A$ is called 
bounded if there is a $\A_\tau^m$ such that $B\subset \A_\tau^m$ and 
$\sup_{a\in B}\{|a|_l\} < \infty$ for any $l=0,1,\cdots$.

\begin{defi}
For any $a(\xi,x)\in\A$, we define the oscillatory integral
$\os[e^{-i\xi x}a]$ by 
\begin{align*}
\os[e^{-i\xi x}a]
&\equiv \osint e^{-i\xi x}a(\xi,x)d\xi dx \\
&=\frac{1}{(2\pi)^n}\lim_{\varepsilon\to 0}\int \chi(\varepsilon\xi,\varepsilon x) 
e^{-i\xi x}a(\xi,x)d\xi dx,
\end{align*}
where $\chi(\xi,x)$ is any function of $\CS(\R_{\xi,x}^{2n})$ satisfying 
$\chi(0,0)=1$.
\end{defi}

\begin{lem}\label{lem:chi}
If $\chi(x)\in\CS(\R^n)$ satisfies $\chi(0)=1$ then
\begin{align}
\chi(\varepsilon x)&\mathop{\to}_{\varepsilon \to 0} 1 
 \quad \text{(uniformly on compact sets)}, \\
\partial^I_x\chi(\varepsilon x)&\mathop{\to}_{\varepsilon \to 0} 0 
 \quad \text{(uniformly on $\R^n$, $|I|>0$)},\label{eq:chiConv}
\end{align}
and for any multi-index $I$ there is a constant $C_I$ independent of 
\ $0<\varepsilon<1$ such that for any $\sigma,0\leq\sigma\leq|I|$ 
\begin{align}
|\partial^I_x \chi(\varepsilon x)|\leq 
C_I\varepsilon^\sigma \langle x \rangle^{-(|I|-\sigma)}.
\label{in:chiEst}
\end{align}
\end{lem}
\begin{proof}
(\ref{eq:chiConv}) is clear since $\partial_x^I \chi(\varepsilon x)
=\varepsilon^{|I|}\partial_y^I \chi(y)|_{y=\varepsilon x}$. 
If $|x|\leq 1$ then (\ref{in:chiEst}) is obtained from the equality 
$|\partial_x^I \chi(\varepsilon x)| = 
\varepsilon^{\sigma}(\varepsilon^{(|I|-\sigma)}|\partial_y^I \chi(y)|_{|y=\varepsilon x})$. 
If $|x|>1$, 
\begin{align*}
\varepsilon^{(|I|-\sigma)}|\partial_y^I \chi(y)_{|y=\varepsilon x}|
 &=(|y|^{(|I|-\sigma)}|\partial_y^I\chi(y)|)_{|y=\varepsilon x}|x|^{-(|I|-\sigma)}\\ 
 &\leq C_I\langle x \rangle ^{-(|I|-\sigma)}\quad (0\leq \sigma \leq |I|).
\end{align*}
\end{proof}

\begin{thm}
For any $a\in\A$, $\os[e^{-i\xi x}a]$ is 
independent of the choice of $\chi \in \CS$ satisfying $\chi(0,0)=1$. 
For $a\in \A_\tau^m$, if we take integers $l,l'$ satisfying
\begin{align}
-2l+m<-n,\quad -2l'+\tau<-n,\label{cond:l1}
\end{align}
then 
\begin{align*}
|\langle x \rangle^{-2l'}\langle D_\xi \rangle^{2l'}
 \{\langle \xi \rangle^{-2l}\langle D_x \rangle^{2l}a(\xi,x)\}|
 \in L_1(\R^{2n}),
\end{align*}
and 
\begin{align*}
\os[e^{-i\xi x}a]=\osint e^{-i\xi x}\langle x \rangle^{-2l'}\langle D_\xi \rangle^{2l'}
 \{\langle \xi \rangle^{-2l}\langle D_x \rangle^{2l}a(\xi,x)\} d\xi dx.
\end{align*}
Moreover, for $a\in \A_\tau^m$ there is a constant $C$ such that
\begin{align}
|\os[e^{-i\xi x}a]|\leq C|a|_{2(l+l')}.\label{cond:bound}
\end{align}
\end{thm}

\begin{proof}
Fix $0<\varepsilon<1$. Performing integrations by parts, we have 
\begin{align*}
I_\varepsilon 
&\equiv \int e^{-i \xi x}\chi(\varepsilon \xi,\varepsilon x)a(\xi,x)d\xi dx \\
&=\int e^{-i \xi x}\langle\xi\rangle^{-2l}\langle D_x\rangle^{2l}(\chi(\varepsilon \xi,\varepsilon x)a(\xi,x))d\xi dx \\
&=\int e^{-i \xi x}\langle x \rangle^{-2l'}\langle D_\xi\rangle^{2l'}
\{ \langle\xi\rangle^{-2l}\langle D_x\rangle^{2l}
(\chi(\varepsilon \xi,\varepsilon x)a(\xi,x))\}d\xi dx. 
\end{align*}
Lemma \ref{lem:chi} implies the set 
$\{\chi(\varepsilon\xi,\varepsilon x) \}_{0<\varepsilon <1}$ is a bounded subset of 
$\A_0^0$, so that for any $I$ and $J$ there is a constant $C_{I,J}$ 
independent of $\varepsilon, a\in \A_\tau^m$ such that
\begin{align*}
|\partial_\xi^I\partial_x^J(\chi(\varepsilon\xi,\varepsilon x)a(\xi,x)|
\leq C_{I,J}|a|_{(|I|+|J|)}\langle\xi\rangle^m\langle x \rangle^\tau.
\end{align*}
On the other hand, for any $s$, there is a constant $C_{s,I}$ such that
\begin{align*}
|\partial_\xi^I\langle\xi\rangle^s|\leq C_{s,I}\langle\xi\rangle^{s-|I|},
\end{align*}
that is obtained by induction from 
$\partial_{\xi_j}\langle \xi \rangle^s=s\xi_j\langle\xi\rangle^{s-2}$. 
From these facts, for any $I$, there is a constant $C_{l,I}$ independent of 
$\varepsilon,a\in\A_\tau^m$ such that
\begin{align*}
|\partial_\xi^I\{ \langle\xi\rangle^{-2l}\langle D_x\rangle^{2l}
(\chi(\varepsilon \xi,\varepsilon x)a(\xi,x))\}|
\leq C_{l,I}|a|_{(2l+|I|)}\langle\xi\rangle^{m-2l}\langle x\rangle^\tau.
\end{align*}
Hence there is a constant $C_{l,l'}$ independent of $\varepsilon,a\in\A_\tau^m$
such that
\begin{align}
\langle x \rangle^{-2l'}\langle D_\xi\rangle^{2l'}
\{ \langle\xi\rangle^{-2l}\langle D_x\rangle^{2l}
(\chi(\varepsilon \xi,\varepsilon x)a(\xi,x))\}
\leq C_{l,l'}|a|_{(l+l')}\langle\xi\rangle^{m-2l}\langle x \rangle^{\tau-2l'}.
\label{ineq:int}
\end{align}
The right-hand side of (\ref{ineq:int}) is in $L_1(\R_{\xi,x}^{2n})$ because 
of Condition (\ref{cond:l1}). Hence Lebesgue's convergence theorem gives
\begin{align*}
\os[e^{i\xi x a}]&=\lim_{\varepsilon \to 0} \frac{I_\varepsilon}{(2\pi)^2}\\
&=\frac{1}{(2\pi)^2}
\int e^{-i \xi x}\langle x \rangle^{-2l'}\langle D_\xi\rangle^{2l'}
\{ \langle\xi\rangle^{-2l}\langle D_x\rangle^{2l}a(\xi,x))\}d\xi dx,
\end{align*}
and proves (\ref{cond:bound}).
\end{proof}

\begin{thm}
Assume $\{a_j\}_{j=1}^\infty$ is a bounded set of $\A$ and there is $a\in \A$ 
such that
\begin{align*}
a_j(\xi,x) \to a(\xi,x) \quad\text{uniformly on compact sets of }\R^{2n}.
\end{align*}
Then 
\begin{align*}
\lim_{j\to\infty} \os[e^{-i\xi x}a_j]=\os[e^{-i\xi x}a]
\end{align*}
holds.
\end{thm}

\begin{thm}
The oscillatory integral satisfies the following formula
\begin{align*}
\os[e^{-i\xi x}a(\xi,x)]&=\os[e^{-i(\xi-\xi_0) (x-x_0)}a(\xi-\xi_0,x-x_0)], \quad (\xi_0,x_0)\in \R^{2n} \\
\os[e^{-i\xi x}x^I a]&=\os[(-D_\xi)^I e^{-i\xi x}a]=\os[e^{-i\xi x}D_\xi^I a],\\
\os[e^{-i\xi x}\xi^I a]&=\os[(-D_x)^I e^{-i\xi x}a]=\os[e^{-i\xi x}D_x^I a].
\end{align*}
\end{thm}

\begin{thm}
Let $a=a(x)\in \A$ be a function depending only on $x$. Then 
\begin{align*}
\osint e^{-i\xi (x-y)} a(x) d\xi dx = a(y).
\end{align*}
\end{thm}

\section{ The Fedosov construction of star products }
We provide here a brief summary of the Fedosov construction that is 
one of the most useful method of constructing a star product on a 
symplectic manifold $(M,\omega)$. 
For details see \cite{Fedosov:1994,Fedosov:1996}.

A formal Weyl algebra $W_x$ associated with $T_x M$ for $x\in M$ is an 
associative algebra with unit over $\C$ defined as follows: Each element of 
$W_x$ is a formal power series in $\lambda$ with coefficients being formal 
polynomials in $T_x M$, that is, each element has the form 
\begin{align*}
a(y,\lambda)=\sum_{k,J}\lambda^k a_{k,J} y^J,
\end{align*}
where $y=(y^1,\ldots,y^{2n})$ are linear coordinates on $T_x M$, 
$J=(j_1,\ldots,j_{2n})$ is a multi-index and 
$y^J=(y^1)^{j_1}\cdots (y^{2n})^{j_{2n}}$. The product $\circ$ is defined 
by the Moyal-Weyl rule,
\begin{align*}
a \circ b = 
\sum_{k=0}^{\infty }\left( \frac{\lambda }{2}\right)
^{k}\frac{1}{k!}\omega ^{i_{1}j_{1}}\cdots \omega ^{i_{k}j_{k}}\frac{
\partial ^{k}a}{\partial y^{i_{1}}\cdots \partial y^{i_{k}}}\frac{\partial
^{k}b}{\partial y^{j_{1}}\cdots \partial y^{j_{k}}},
\end{align*}
where $\omega^{lm}$ are the coefficients of $\omega$ with respect to $y^j$. 
If we assign $\deg(y^j)=1$ and $\deg(\lambda)=2$, the algebra 
$W_x$ becomes a filtered algebra.

Let $W=\cup _{x\in M}W_{x}$. Then $W$ is a bundle of algebras over $M$,
called the Weyl bundle over $M$. Each section of $W$ has the form
\begin{align}  \label{eq:weylSection}
a(x,y,\lambda )=\sum_{k,\alpha }\lambda ^{k}a_{k,\alpha }(x)y^{\alpha },
\end{align}
where $x\in M$. 
We call $a(x,y,\lambda )$ smooth if each coefficient 
$a_{k,\alpha }(x)$ is smooth in $x$. 
We denote the set of smooth sections by $\Gamma W$. 
It constitutes an associative algebra with unit under the fibrewise
multiplication.

Let $\nabla$ be a torsion-free symplectic connection on $M$, 
which always exists 
and $\partial:\Gamma W\to\Gamma W\otimes\Lambda^1$ be its induced 
covariant derivative. 
Consider a connection on $W$ of the form
\begin{align}  \label{eq:weylConnection}
Da=-\delta a+\partial a-\frac{1}{\lambda }[\gamma ,a],\quad 
\text{ for }a\in\Gamma W
\end{align}
with $\gamma \in \Gamma W\otimes \Lambda ^{1}$, where
\begin{align*}
\delta a=dx^{k}\wedge \frac{\partial a}{\partial y^{k}}.
\end{align*}
Clearly, $D$ is a derivation for the Moyal-Weyl product.
A simple computation shows that
\begin{align*}
D^2a=\frac{1}{\lambda}[\Omega,a],\quad\text{for any }a\in \Gamma W,
\intertext{where}
\Omega=\omega-R+\delta\gamma-\partial\gamma+\frac{1}{\lambda}\gamma^2.
\end{align*}
Here $R=\frac{i}{4}R_{ijkl}y^{i}y^{j}dx^{k}\wedge dx^{l}$ and 
$R_{ijkl}=\omega_{im}R_{jkl}^{m}$ is the curvature tensor of the 
symplectic connection.

A connection of the form (\ref{eq:weylConnection}) is called 
Abelian if $\Omega$ is a scalar 2-form, that is, 
$\Omega\in\Lambda^2[[\lambda]]$. 
We call $D$ a Fedosov connection if it is Abelian and 
$\deg \ \gamma\geq 3$. 
For an Abelian connection, the Bianchi identity implies that 
$d\Omega=D\Omega=0$, that is, $\Omega$ is closed. 
In this case, we call $\Omega$ a Weyl curvature.

\begin{thm}
[\protect\cite{Fedosov:1994}]
\label{thm:fedosovConnection} Let $\nabla $ be any
torsion-free symplectic connection, and 
$\Omega =\omega +\lambda \omega_{1}+\cdots \in Z^{2}(M)[[\lambda ]]$ 
a perturbation of the symplectic form
$\omega$.
There exists a unique $\gamma \in \Gamma W\otimes \Lambda ^{1}$ 
such that $D$ given by Equation (\ref{eq:weylConnection}) is a 
Fedosov connection which has Weyl curvature $\Omega $ and satisfies 
$\delta ^{-1}\gamma =0$.
\end{thm}

We denote $W_{D}$ the set of smooth and flat sections, that is, 
a section $a\in \Gamma W$ satisfying $Da=0$. 
The space $W_{D}$ becomes a subalgebra of $\Gamma W$. 
Let $\sigma $ denote the projection of 
$W_{D}$ onto $C^\infty(M)[[\lambda ]]$ defined by $\sigma (a)=a|_{y=0}$.

\begin{thm}
[\protect\cite{Fedosov:1994}] 
\label{thm:FedosovQuantization} 
Let $D$ be an Abelian connection. 
Then, for any $a_0(x,\lambda)\in C^\infty(M)[[\lambda]]$ there 
exists a unique section $a\in W_D$such that $\sigma(a)=a_0$. 
Therefore, $\sigma$ establishes an isomorphism between $W_D$ and 
$C^\infty(M)[[\lambda]]$ as $\C [[\lambda]]$-vector spaces.
\end{thm}

We denote the inverse map of $\sigma $ by $Q$ and call it a 
quantization procedure. 
The Weyl product $\circ$ on $W_{D}$ is translated to 
$C^\infty(M)[[\lambda]]$ hence yielding a star product $*$. 
Namely, we set for $a,b\in C^\infty(M)[[\lambda ]]$
\begin{align*}
a*b=\sigma (Q(a) \circ Q(b)).
\end{align*}

For $G$-invariant star products, there is a simple criterion.

\begin{prop}
[\protect\cite{Fedosov:1996}]
\label{prop:GinvariantStar}
Let $\nabla $ be a $G$-invariant connection, $\Omega $ be a $G$-invariant
Weyl curvature and $D$ be the Fedosov connection corresponding to 
$(\nabla,\Omega)$. 
Then the star product corresponding to $D$ is $G$-invariant.
\end{prop}

We study mainly star products of Fedosov type because of the 
following theorem.

\begin{thm}
[\protect\cite{BBG:1998}]
\label{thm:GStarClass}
Every $G$-invariant star product is $G$-equivalent to a Fedosov star product.
\end{thm}

\section*{Acknowledgements}

I would like to thank Giuseppe Dito for many fruitful discussions 
on the subject and Izumi Ojima for many important suggestions. 
It is pleasure to thank Daniel Sternheimer for warmest hospitality. 

This research was supported by a Postdoctoral Scholarship of 
the Ministry of National Education and Research of France.

\bibliographystyle{plain}

\end{document}